\documentclass[12pt]{amsart}
\usepackage{amsmath,amscd,amssymb,amsfonts}
\newcommand{\h}{\hbox}
\newcommand{\q}{\quad}
\newcommand{\nin}{\par\noindent}
\newcommand{\bs}{\par\bigskip}
\newcommand{\ms}{\par\medskip}
\newcommand{\sk}{\par\smallskip}
\newcommand{\bl}{\bigl}
\newcommand{\br}{\bigl}
\newcommand{\ssb}{\raise.15ex\h{${\scriptscriptstyle\bullet}$}}
\newcommand{\scc}{\,\raise.15ex\h{${\scriptstyle\circ}$}\,}
\newcommand{\msum}{\h{$\sum$}}
\newcommand{\mopl}{\h{$\bigoplus$}}
\newcommand{\mcap}{\h{$\bigcap$}}
\newcommand{\mcup}{\h{$\bigcup$}}
\newcommand{\mprod}{\h{$\prod$}}
\newcommand{\al}{\alpha}
\newcommand{\C}{{\mathbf C}}
\newcommand{\e}{{\mathbf e}}
\newcommand{\ep}{\varepsilon}
\newcommand{\E}{\widetilde{E}}
\newcommand{\Ec}{{\mathcal E}}
\newcommand{\Hc}{{\mathcal H}}
\newcommand{\Ht}{\widetilde{H}}
\newcommand{\HH}{{\mathbf H}}
\newcommand{\I}{{\mathcal I}}
\newcommand{\La}{\Lambda}
\newcommand{\la}{\lambda}
\newcommand{\Lc}{{\mathcal L}}
\newcommand{\M}{{\mathcal M}}

\newcommand{\MM}{{\mathcal M}^{\ssb}}
\newcommand{\N}{{\mathbf N}}
\newcommand{\Oc}{{\mathcal O}}
\newcommand{\PP}{{\mathbf P}}
\newcommand{\Q}{{\mathbf Q}}
\newcommand{\Qt}{\widetilde{Q}}
\newcommand{\R}{{\mathbf R}}
\newcommand{\St}{\widetilde{S}}
\newcommand{\So}{{}\,\overline{\!S}{}}
\newcommand{\Sc}{{\mathcal S}}

\newcommand{\Si}{\Sigma}
\newcommand{\si}{\sigma}
\newcommand{\X}{\widetilde{X}}
\newcommand{\Z}{{\mathbf Z}}
\newcommand{\BM}{{\rm BM}}
\newcommand{\CH}{{\rm CH}}
\newcommand{\DR}{{\rm DR}}
\newcommand{\Gr}{{\rm Gr}}
\newcommand{\MHM}{{\rm MHM}}
\newcommand{\Sing}{{\rm Sing}}
\newcommand{\Sp}{{\rm Sp}}
\newcommand{\vir}{{\rm vir}}
\newcommand{\codim}{{\rm codim}}
\newcommand{\into}{\hookrightarrow}
\newcommand{\simto}{\buildrel\sim\over\longrightarrow}
\newcommand{\ges}{\geqslant}
\newcommand{\les}{\leqslant}
\begin{document}
\title[Hirzebruch-Milnor classes and Steenbrink spectra] {Hirzebruch-Milnor classes and Steenbrink spectra\\ of certain projective hypersurfaces}
\author[L. Maxim ]{Laurentiu Maxim}
\address{L. Maxim : Department of Mathematics, University of Wisconsin-Madison, 480 Lincoln Drive, Madison WI 53706-1388 USA}
\email{maxim@math.wisc.edu}
\author[M. Saito ]{Morihiko Saito}
\address{M. Saito: RIMS Kyoto University, Kyoto 606-8502 Japan}
\email{msaito@kurims.kyoto-u.ac.jp}
\author[J. Sch\"urmann ]{J\"org Sch\"urmann}
\address{J. Sch\"urmann : Mathematische Institut, Universit\"at M\"unster, Einsteinstr. 62, 48149 M\"unster, Germany}
\email{jschuerm@uni-muenster.de}
\dedicatory{To the memory of Friedrich Hirzebruch}
\begin{abstract} We show that the Hirzebruch-Milnor class of a projective hypersurface, which gives the difference between the Hirzebruch class and the virtual one, can be calculated by using the Steenbrink spectra of local defining functions of the hypersurface if certain good conditions are satisfied, e.g. in the case of projective hyperplane arrangements, where we can give a more explicit formula. This is a natural continuation of our previous paper on the Hirzebruch-Milnor classes of complete intersections.
\end{abstract}
\maketitle
\centerline{\bf Introduction}
\bs\nin
In his classical book \cite{Hi}, F.~Hirzebruch introduced the cohomology Hirzebruch characteristic class $T_y^*(TX)$ of the tangent bundle $TX$ of a compact complex manifold $X$, see also (1.1) below. It belongs to $\HH^{\ssb}(X)[y]$ where $\HH^k(X)=H^{2k}(X,\Q)$. By specializing to $y=-1$, $0$, $1$, it specializes to the Chern class $c^*(TX)$, the Todd class $td^*(TX)$, and the Thom-Hirzebruch $L$-class $L^*(TX)$ respectively, see \cite[Sect.~5.4]{HBJ}. Its highest degree part, which is called the $T_y$-genus in \cite[10.2]{Hi}, was mainly interested there by the relation with his Riemann-Roch theorem. This coincides with the $\chi_y$-genus
$$\chi_y(X):=\msum_p\,\chi(\Omega_X^p)\,y^p\in\Z[y],$$
which specializes to the Euler characteristic, the arithmetic genus, and the signature of $X$ for $y=-1$, $0$, $1$.
The cohomology class $T_y^*(TX)$ is identified by Poincar\'e duality with the homology Hirzebruch class $T_{y*}(X)$ in the smooth case. It is generalized to the singular case by \cite{BSY} (see (1.2) below).
\sk
Hirzebruch also introduced there the virtual $T_y$-genus (or $T_y$-characteristic) which gives the $T_y$-genus of smooth complete intersections $X$ in smooth projective varieties $Y$. We can define the virtual Hirzebruch class $T^{\,\vir}_{y*}(X)$ of any complete intersection $X$ like the virtual $T_y$-genus even if $X$ is singular. In this paper we adopt a more sophisticated construction as in \cite{BSY}, see (1.2--3) below. This is compatible with the construction in \cite[11.2]{Hi}, see (1.4) below. In \cite{MSS} we proved that the difference between the Hirzebruch class and the virtual one is given by the {\it Hirzebruch-Milnor class} $M_y(X)$ supported on the singular locus of $X$, and gave an inductive formula in the case of global complete intersections with arbitrary singularities.
\sk
In this paper, we restrict to the case of projective hypersurfaces (i.e. the codimension is one) satisfying certain good conditions in order to prove a formula for the Hirzebruch-Milnor class $M_y(X)$ by using the Steenbrink spectra (see \cite{St1}, \cite{St2}) of local defining functions of $X$ in $Y$. This is a natural continuation of the last section of \cite{MSS}. (Note that the implication of the calculations in loc.~cit.\ was not explained there).
More precisely, let $Y$ be a smooth complex projective variety having a very ample line bundle $L$. Set
$$X=s^{-1}(0)\q\h{with}\,\,\,s\in\Gamma(Y,L^{\otimes m})\,\,\,\h{for some positive integer}\,\,\,m.$$
Let $s'_1,\dots,s'_{n+1}$ be sufficiently general sections of $L$, where $n:=\dim X$. Take sufficiently general non-zero complex numbers  $a_j$ with $|a_j|$ sufficiently small ($j\in[1,n])$. For $j\in[1,n+1]$, set
$$\aligned s_{a,j}&:=s-a_1s^{\prime\,m}_1-\cdots-a_{j-1}s^{\prime\,m}_{j-1},\q f_{a,j}:=\bl(s_{a,j}/s_j^{\prime\,m})|_{Y\setminus X'_j},\q X'_j:=s^{\prime\,-1}_j(0),\\
X_{a,j}&:=s_{a,j}^{-1}(0),\q\Si_j:=\Sing\,X_{a,j}\,(=\mcap_{k<j}X'_k\cap\Si),\q \Si:=\Si_1.\endaligned$$
Set $r:=\max\{j\mid\Si_j\ne\emptyset\}$.
By \cite{MSS}, \cite{MSS2}, there is the Hirzebruch-Milnor class $M_y(X)\in\HH_{\ssb}(\Si)[y]$ with $\HH_k(\Si)=H_{2k}^{\BM}(\Si,\Q)$ or $\CH_k(\Si)_{\Q}$, satisfying
$$T^{\,\vir}_{y*}(X)-T_{y*}(X)=(i_{\Si,X})_*M_y(X),
\leqno(0.1)$$
\vskip-20pt
$$M_y(X)=\msum_{j=1}^r\,T_{y*}\bl((i_{\Si_j\setminus X'_j,\Si})_{!\,}\varphi_{f_{a,j}}\Q_{h,Y\setminus X'_j}\br),
\leqno(0.2)$$
where $i_{A,B}:A\into B$ denotes the inclusion for $A\subset B$ in general.
In this paper, $\varphi_{f_{a,j}}\Q_{h,Y\setminus X'_j}$ denotes a mixed Hodge module up to a shift of complex on $\Si_j\setminus X'_j$ such that its underlying $\Q$-complex is the vanishing cycles $\varphi_{f_{a,j}}\Q_{Y\setminus X'_j}$ in \cite{De3}, see \cite{Sa1}, \cite{Sa2}. For the definition of $T_{y*}(\MM)$ with $\MM$ a bounded complex of mixed Hodge modules, see (1.2.1) below.
In \cite{MSS} we assumed $m=1$, but it is not difficult to generalize the argument there to the case $m>1$ by using (2.4) below, see \cite{MSS2} for details.
By \cite[Prop.~5.21]{Sch}, this formula specializes at $y=-1$ to a formula for the Chern classes, which was conjectured by S.~Yokura \cite{Yo2}, and was proved by A.~Parusi\'nski and P.~Pragacz \cite{PP} (where $m=1$).
\sk
In the case of hyperplane arrangements, we may assume
$$f_{a,j}=g_j+x_1^m+\cdots+x_{j-1}^m,$$
for some coordinates $x_1,\dots,x_{n+1}$ of $Y\setminus X'_j=\C^{n+1}$, where $g_j$ defines the restriction of the hyperplane arrangement to $Y\setminus X'_j$. We have a topologically trivial one-parameter family
$$g_j(\la x_1,\dots,\la x_{j-1},x_j,\dots,x_{n+1})+x_1^m+\cdots+x_{j-1}^m\q(\la\in\C),$$
and apply the Thom-Sebastiani theorem \cite{Sa3} at $\la=0$ (together with \cite{DMST}).
This argument can be extended to the general case by using the deformation to the normal cone. Set $Y_j:=\mcup_{k<j}\,X'_k$, $g_{a,j}:=f_{a,j}|_{Y_j\setminus X'_j}$. For the calculation of the right-hand side of (0.2), it is then sufficient to calculate $\varphi_{g_{a,j}}\Q_{h,Y_j\setminus X'_j}$ together with the action of the semisimple part of the monodromy $T_s$, see \cite{MSS2} for details.
From now on, we fix $j\in[1,r]$, and denote $g_{a,j}$, $Y_j$, $Y_j\cap X'_j$, $\Si_j$ respectively by $f$, $Y$, $X'$, $\Si$ to simplify the notation.
\sk
Let $\Sc$ be a complex algebraic stratification of $\Si\setminus X'$ such that the $\Hc^j_S:=\Hc^j\varphi_f\Q_{Y\setminus X'}|_S$ are local systems for any strata $S\in\Sc$ (which are assumed smooth). These local systems canonically underlie admissible variations of mixed Hodge structure $\HH^j_S$, since $\varphi_f\C_{Y\setminus X'}$ underlies a mixed Hodge module up to a shift of complex $\varphi_{h,f}\Q_{Y\setminus X'}$. Let $\Hc^j_{S,\la}\subset\Hc^j_S$ be the $\la$-eigenspace by the action of the semisimple part $T_s$ of the Milnor monodromy $T$ (which is defined as the monodromy of the local system on a punctured disk associated with the Milnor fibration, see \cite{De3}).
The local system monodromy of $\Hc^j_{S,x}$ around $X'$ coincides with the $m$-th power of the Milnor monodromy where we take $x\in S$ sufficiently near $X'$ so that we have a loop around $X'$ passing through $x$. (This can be reduced to Lemma~(3.3) below by using the expression $f=(s/s^{\prime\,m})|_{Y\setminus X'}$ together with the deformation to the normal cone of $X'$.) So the situation is quite different from the one in \cite{CMSS} where $X$ is a fiber of a morphism to a curve, and a formula for the Hirzebruch-Milnor class is given by using the mixed Hodge structure on each stalk of the vanishing cycles under the assumption on the triviality of the global monodromies of the local systems.
\sk
In this paper we prove a variant of it by using the Steenbrink spectrum (see (1.5) below):
$$\Sp(f,x)=\msum_{\al\in\Q}\,n_{f,x,\al}t^{\al}\q(x\in\Si\setminus X').$$
Here $n_{f,x,\al}\in\Z$ is independent of $x\in S$, and will be denoted by $n_{f,S,\al}$.
We say that a compactification $\St$ of $S$ is {\it good} if $D_{\St}:=\St\setminus S$ is a divisor with simple normal crossings on a smooth projective variety $\St$ and the natural inclusion $S\into\Si$ extends to a (unique) morphism $\pi_{\St}:\St\to\Si$.
Using the results in the last section of \cite{MSS}, we get the following.
\ms\nin
{\bf Theorem~1.} {\it Assume the following two conditions$\,:$
\sk\nin
$(a)$ Every $\HH^j_S$ is a locally constant variation of mixed Hodge structure on $S$.
\sk\nin
$(b)$ Each $\Hc^j_{S,\la}$ is isomorphic to a direct sum of copies of a rank $1$ local system $L_{S,\la}$ which is independent of $j$.
\sk\nin
Let $\Lc_{\St,\la}$ be the Deligne extension of $L_{S,\la}$ as an $\Oc_{\St}$-module with a logarithmic connection such that the eigenvalues of the residues of the connection are contained in $(0,1]$ where $\St$ is any good compactification of $S$. Then
$$\aligned&\q T_{y*}\bl((i_{\Si\setminus X',\Si})_{!\,}\varphi_f\Q_{h,Y\setminus X'}\br)\\
&=\sum_{S,\al,q}(-1)^{q+n-1}\,n_{f,S,\al}\,(\pi_{\St})_*\,td_{(1+y)*}\bl[\Lc_{\St,\e(-\al)}\otimes_{\Oc_{\St}}\Omega^q_{\St}(\log D_{\St})\br](-y)^{\lfloor n-\al\rfloor+q},\endaligned
\leqno(0.3)$$
where $\e(\al):=\exp(2\pi i\al)$, $\lfloor n-\al\rfloor$ denotes the integer part $($see $(1.5.2)$ below$)$, and $td_{(1+y)*}$ is as in $(1.2.2)$ below.}
\ms
Here the sign $(-1)^{n-1}$ comes from the definition of spectrum, see (1.5.1) below.
Note that the assertion for the Chern-Milnor class $M(X)$ corresponding to Theorem~1 (or deduced from it by specializing to $y=-1$) is essentially a corollary of \cite{PP}, and holds without assuming conditions~$(a)$, $(b)$, see (2.5) below.
\sk
Condition~$(a)$ in Theorem~1 is satisfied if $X$ is locally analytically trivial along each stratum $S$ (e.g. if the intersection of $X$ with transversal slices to any $S$ has only isolated singularities of type $A$, $D$, $E$) or if the Hodge filtration $F$ on any $\Hc^j_{S,\la}\otimes_{\C}\Oc_S$ is trivial (e.g. if every nonzero $\Hc^j_{S,\la}$ has rank 1).
As for condition~$(b)$, we have the following.
\ms\nin
{\bf Proposition~1.} {\it Condition~$(b)$ is satisfied if the following two conditions hold$\,:$
\sk\nin
$(c)$ Every $S\in\Sc$ has a simply connected good compactification $\St$.
\sk\nin
$(d)$ The local monodromy of $\Hc^j_{S,\la}$ around each irreducible component $D_{{\St},i}$ of $D_{\St}=\St\setminus S$ is the multiplication by a constant number $c_{S,i,\la}$ which is independent of $j$.}
\ms
Conditions~$(a)$, $(c)$, $(d)$ are satisfied, for instance, in the case $X$ is a projective hyperplane arrangement in $\PP^n$ (see Proposition~3 below), or the projective compactification of the affine cone of a hypersurface in $\PP^{n-1}$ which has only isolated singularities with semisimple Milnor monodromies. In these examples, the following conditions are satisfied for any $S\in\Sc:$
$$c_{S,i,\la}=\la^{-m'_{S,i}}\q\h{with}\q m'_{S,i}\in\Z,
\leqno(0.4)$$
\vskip-6mm
$$X\cap Z_S=g_S^{-1}(0)\,\,\,\h{with $\,g_S\,$ a homogeneous polynomial,}
\leqno(0.5)$$
where $Z_S\subset Y$ is an analytic transversal slice to $S$ which intersects $S$ transversally at a sufficiently general point of $S$. We assume that $Z_S$ is sufficiently small, and has some coordinates to express $g_S$.
\sk
Set
$$\La_S:=\{\la\in\C^*\mid \Hc^j_{S,\la}\ne 0\,\,(\exists\,j\in\N)\}.$$
Let $G_S\subset\C^*$ be the subgroup generated by $\la\in\La_S$. It corresponds to a finitely generated $\Z$-submodule of $\Q$ by $\al\mapsto\e(\al)=\exp(2\pi i\al)$, and is generated by $\e(1/m'_S)$ with $m'_S\in\Z_{>0}$. Note that the $m'_{S,i}$ are well-defined mod $m'_S$, and we sometimes assume
$$m'_{S,i}/m'_S\in[0,1).
\leqno(0.6)$$
However, it is not necessarily easy to give $m'_S$ {\it explicitly} in general (even in the hyperplane arrangement case).
If condition~(0.5) is satisfied, then $\Hc^j_{S,\la}=0$ unless $\la^{m_S}=1$, and we get
$$m_S:=\deg g_S\in\Z\,m'_S.
\leqno(0.7)$$
Here the equality $m_S=m'_S$ does not always hold (e.g. if $f=y_1^2y_2^2$ with $m_S=4$, $m'_S=2$). We will assume (0.6) with $m'_S$ replaced by $m_S$ in case $m'_S$ is not explicitly given.
\ms\nin
{\bf Proposition~2.} (i) {\it Assume conditions~$(c)$, $(d)$ and $(0.4)$, $(0.6)$ hold. Then there is a rank $1$ local system $L'_S$ on $S$ such that the eigenvalues of its local monodromies $c_{S,i,\la}$ satisfy $(0.4)$ with $\la=\e(1/m'_S)$ and we have $L_{S,\la}\cong L_S^{\prime\,\otimes k}$ for any $\la=\e(k/m'_S)\in\La_S$.
Let $\Lc'_{\St}$ be the Deligne extension of $L'_S$ on $\St$ such that the eigenvalues of the residues are contained in $[0,1)$. Then
$$\Lc_{\St,\la}=\Lc_{\St}^{\prime\,\otimes k}\otimes_{\Oc_{\St}}\Oc_{\St}\bl(\msum_i\,(\lceil k\,m'_{S,i}/m'_S\rceil-1)D_{{\St},i}\br)\q\h{for}\,\,\,\la=\e(k/m'_S)\in\La_S,
\leqno(0.8)$$
where $\lceil*\rceil$ is as in $(1.5.2)$ below.}
\sk\nin
(ii) {\it With the above assumption, assume further that condition~$(0.5)$ holds and there is a rank $1$ local system $L_S$ on $S$ such that the eigenvalues of its local monodromies $c_{S,i,\la}$ satisfy $(0.4)$ with $\la=\e(1/m_S)$. Then the assertion of {\rm (i)} together with $(0.6)$ holds with $L'_S$, $\Lc'_{\St}$, $m'_S$ replaced respectively by $L_S$, $\Lc_{\St}$, $m_S$.}
\ms
There is a certain similarity between (0.8) and \cite[1.4.3]{BS2}. The difference between $\lceil*\rceil-1$ in (0.8) and $\lfloor*\rfloor$ in loc.~cit.\ comes from the difference between $j_!$ and $\R j_*$ if $j$ denotes the open embedding $S\into\St$. (It is also related to \cite[Thm.~4.2]{BS1} and \cite[3.10.9]{Sa2}).
\sk
In the hyperplane arrangement case, an explicit formula for the $n_{f,S,\al}$ is given by \cite{BS2} in the case $X$ is reduced and $\codim_Y\,S\les 3$, and we have the following (see Propositions~(3.2) and (3.7) below).
\ms\nin
{\bf Proposition~3.} {\it If $X$ is a projective hyperplane arrangement in $\PP^n$, then conditions~$(a)$, $(c)$, $(d)$ and $(0.4)$, $(0.5)$ are all satisfied, and $L_S$ exists so that $(0.6)$ and $(0.8)$ hold with $\Lc'_{\St}$, $m'_S$ replaced by $\Lc_{\St}$, $m_S$. Moreover, $\St$, $\Lc_{\St}$, $m'_{S,i}$ are described explicitly, and the Hirzebruch-Milnor class is an combinatorial invariant of the hyperplane arrangement.}
\ms
The proof of the last assertion follows from an argument similar to \cite{BS2}, where the combinatorial property of the spectrum is shown by using the Hirzebruch-Riemann-Roch theorem together with \cite{DP}.
Here we use $\HH_k(X):=\CH_k(X)_{\Q}$, since the structure of $\CH_k(\Si)_{\Q}$ is quite simple (see Proposition~(3.6) below).
\ms
The first named author is partially supported by NSF-1005338. The second named author is partially supported by Kakenhi 24540039. The third named author is supported by the SFB 878 ``groups, geometry and actions''.
\ms
In Section 1 we review some basics of Hirzebruch characteristic classes and Steenbrink spectra of hypersurfaces. In Section 2 we give the proofs of Theorem~1 and Propositions~1 and 2. In Section 3 we treat the hyperplane arrangement case.
\bs\bs
\vbox{\centerline{\bf 1. Preliminaries}
\bs\nin
In this section we review some basics of Hirzebruch characteristic classes and Steenbrink spectra of hypersurfaces.}
\ms\nin
{\bf 1.1.~Cohomology Hirzebruch classes.} In \cite{Hi}, Hirzebruch introduced the cohomology Hirzebruch characteristic class $T_y^*(TX)$ of the tangent bundle $TX$ of a compact complex manifold $X$ of dimension $n$. By using the formal Chern roots $\{\al_i\}$ for $TX$ satisfying
$$\mprod_{i=1}^n(1+\al_it)=\msum_{j=0}^n\,c_j(TX)\,t^j,$$
it can be defined by
$$T^*_y(TX):=\mprod_{i=1}^n\,Q_y(\al_i)\in\HH^{\ssb}(X)[y].
\leqno(1.1.1)$$
Here we use the formal power series in $\Q[y][[\al]]$
$$Q_y(\al):=\frac{\al(1+y)}{1-e^{-\al(1+y)}}-\al y,\q\Qt_y(\al):=\frac{\al(1+ye^{-\al})}{1-e^{-\al}},\q R_y(\al)=\frac{e^{\al(1+y)}-1}{e^{\al(1+y)}+y},
\leqno(1.1.2)$$
see \cite[10.2 and 11.1]{Hi}. These have the relation
$$Q_y(\al)=(1+y)^{-1}\,\Qt_y(\al(1+y))=\al/R_y(\al).
\leqno(1.1.3)$$
\sk
By specializing to $y=-1$, $0$, $1$, the power series $Q_y(\al)$ becomes respectively
$$1+\al,\q\al/(1 -e^{-\al}),\q\al/\tanh\al,$$
and hence $T^*_y(TX)$ specializes to the Chern class $c^*(TX)$, the Todd class $td^*(TX)$, and the Thom-Hirzebruch $L$-class $L^*(TX)$, see \cite[Sect.~5.4]{HBJ}.
\ms\nin
{\bf 1.2.~Homology Hirzebruch classes.} The cohomology class $T_y^*(TX)$ is identified by Poincar\'e duality with the (Borel-Moore) homology class $T_y^*(TX)\cap[X]$, and this gives the definition of the homology Hirzebruch class $T_{y*}(X)$ in the smooth case. It is generalized to the singular case by \cite{BSY}. Here we can use either the du Bois complex in \cite{dB} or the bounded complex of mixed Hodge modules $\Q_{h,X}$ whose underlying $\Q$-complex is the constant sheaf $\Q_X$ in \cite{Sa2}.
\sk
Let $\MHM(X)$ be the abelian category of mixed Hodge modules on a complex algebraic variety $X$ (see \cite{Sa1}, \cite{Sa2}). For $\MM\in D^b\MHM(X)$, its homology Hirzebruch characteristic class is defined by
$$\aligned T_{y*}(\MM)&:=td_{(1+y)*}\bl(\DR_y[\MM]\br)\in\HH_{\ssb}(X)\bl[y,\h{$\frac{1}{y(y+1)}$}\br]\q\h{with}\\
\DR_y[\MM]&:=\msum_{i,p}\,(-1)^i\,\bl[\Hc^i\Gr_F^p\DR(\MM)\br]\,(-y)^p\in K_0(X)[y,y^{-1}],\raise12pt\h{ }\endaligned
\leqno(1.2.1)$$
where $\HH_k(X):=H^{\BM}_{2k}(X,\Q)$ or $\CH_k(X)_{\Q}$, and
$$td_{(1+y)*}:K_0(X)[y,y^{-1}]\to\HH_{\ssb}(X)\bl[y,\h{$\frac{1}{y(y+1)}$}\br]
\leqno(1.2.2)$$
is given by the scalar extension of the Todd class transformation
$$td_*:K_0(X)\to\HH_{\ssb}(X),$$
(which is denoted by $\tau$ in \cite{BFM}) followed by the multiplication by $(1+y)^{-k}$ on the degree $k$ part (see \cite{BSY}). The last multiplication is closely related with the first equality of (1.1.3). By \cite[Prop.~5.21]{Sch}, we have
$$T_{y*}(\MM)\in\HH_{\ssb}(X)[y,y^{-1}].$$
\sk
The homology Hirzebruch characteristic class $T_{y*}(X)$ of a complex algebraic variety $X$ is defined by applying the above definition to the case $\MM=\Q_{h,X}$ (see \cite{BSY}), i.e.
$$\aligned &T_{y*}(X):=T_{y*}(\Q_{h,X})=td_{(1+y)*}\DR_y[X]\in\HH_{\ssb}(X)[y],\\ &\h{with}\q\DR_y[X]:=\DR_y[\Q_{h,X}].\endaligned$$
This coincides with the definition using the du Bois complex \cite{dB}. It is known that $T_{y*}(X)$ belongs to $\HH_{\ssb}(X)[y]$, see \cite{BSY}. In case $X$ is smooth, we have
$$\DR_y[X]=\La_y[T^*X],
\leqno(1.2.3)$$
where we set for a vector bundle $V$ on $X$
$$\La_y[V]:=\msum_{p\ges 0}\,[\La^pV]\,y^p\in K^0(X)[y],
\leqno(1.2.4)$$
In fact, we have
$$\DR(\Q_{h,X})=\DR(\Oc_X)[-n]=\Omega_X^{\ssb}\q\h{with}\,\,\, n:=\dim X,
\leqno(1.2.5)$$
where the Hodge filtration $F^p$ on $\Omega_X^{\ssb}$ is defined by the truncation $\si_{\ges p}$ as in \cite{De2}.
(For the proof of the coincidence with the above definition of $T_{y*}(X)$ in the smooth case, we have to use the first equality of (1.1.3) and some calculation about Hirzebruch's power series $Q_y(\al)$ as in \cite[Sect.~5.4]{HBJ}, or in the proof of \cite{Yo1}, Lemma~2.3.7, which is closely related with the generalized Hirzebruch-Riemann-Roch theorem as in \cite[Thm.~21.3.1]{Hi}.)
\ms\nin
{\bf 1.3.~Virtual Hirzebruch classes.} Hirzebruch \cite{Hi} also introduced the virtual $T_y$-genus (or $T_y$-characteristic) which gives the $T_y$-genus of smooth complete intersections in smooth projective varieties. Let $X$ be any complete intersection in a smooth projective variety $Y$. We can define the virtual Hirzebruch characteristic class $T^{\,\vir}_{y*}(X)$ by
$$T^{\,\vir}_{y*}(X):=td_{(1+y)*}\DR^{\vir}_y[X]\in\HH_{\ssb}(X)[y],
\leqno(1.3.1)$$
with $\DR^{\vir}_y[X]$ the image in $K_0(X)[[y]]$ of
$$\La_y(T^*_{\vir}X)=\La_y[T^*Y|_X]/\La_y[N^*_{X/Y}]\in K^0(X)[[y]].
\leqno(1.3.2)$$
Here $N^*_{X/Y}$ is defined by the locally free sheaf $\I_X/\I^2_X$ on $X$
with $\I_X\subset\Oc_Y$ the ideal sheaf of the subvariety $X$ of $Y$, and we set for a virtual vector bundle $V$ on $X$
$$\La_yV:=\msum_{p\ges 0}\,\La^pV\,y^p\in K^0(X)[[y]].
\leqno(1.3.3)$$
Note that $\DR^{\vir}_y[X]$ belongs to $K_0(X)[y]$, see \cite[Prop.~3.4]{MSS}.
(We denote by $K^0(X)$ and $K_0(X)$ the Grothendieck group of locally free sheaves of finite length and that of coherent sheaves respectively.)
\sk
We have the equality $T_{y*}(X)=T^{\,\vir}_{y*}(X)$ if $X$ is smooth. So the problem is how to describe their difference in the singular case, and this is given by the Hirzebruch-Milnor class as is explained in the introduction where only the hypersurface case is treated, see \cite{MSS} for the complete intersection case. (For the degree-zero part, i.e.\ on the level of Hodge polynomials, see also \cite{LM}.)
\ms\nin
{\bf 1.4.~Relation with Hirzebruch's construction \cite{Hi}.} 
The image of the above virtual $T_y$-characteristic class of $X$ by the trace morphism ${\rm Tr}_X:\HH_{\ssb}(X)\to\Q$ coincides with the virtual $T_y$-genus of $X$ constructed in \cite[11.2]{Hi}, where $X$ is a (global) complete intersection of codimension $r$ in a smooth complex projective variety $Y$ with $i:X\into Y$ the natural inclusion.
For this we have to recall the cohomological transformation $T_y^*$ (as in \cite{CMSS} in the hypersurface case) applied to the virtual tangent bundle
$$T_{\vir}X:=[TY|_X]-[N_{X/Y}]\in K^0(X).$$
This can be defined by
$$T_y^*(T_{\vir}X):=\frac{\mprod_i\,Q_y(\al_i)}{\mprod_j\,Q_y(\beta_j)}\bigg|_X=\frac{T^*_y(TY)}{\mprod_j\,Q_y(\beta_j)}\bigg|_X.
\leqno(1.4.1)$$
Here the $\al_i$ are the formal Chern roots of $TY$, and the $\beta_j$ ($j\in[1,r]$) are the cohomology classes of hypersurfaces of $Y$ whose intersection is $X$.
\sk
By \cite[Prop.~1.3.1]{MSS} we have the equality
$$T_{y*}^{\vir}(X)=T_y^*(T_{\vir}X)\cap[X].
\leqno(1.4.2)$$
This can be shown by using the first equality of (1.1.3) together with an argument similar to \cite[Lemma~2.3.7]{Yo1}.
By the projection formula, we then get
$$i_*\bigl(T_{y*}^{\vir}(X)\bigr)=\frac{T^*_y(TY)}{\mprod_j\,Q_y(\beta_j)}\cap i_*[X].
\leqno(1.4.3)$$
We have moreover
$$i_*[X]=(\beta_1\cup\dots\cup\beta_r)\cap[Y].
\leqno(1.4.4)$$
This follows, for instance, from the compatibility of the cycle class map $\CH^{\ssb}(Y)\to H^{2\ssb}(Y)$ with the multiplicative structures, see \cite{Fu}.
(Here $X$ is defined scheme-theoretically by using a regular sequence, and we have $[X]=\sum_km_k[X_k]$ with $X_k$ the reduced irreducible components of $X$ and $m_k$ the multiplicities. The equality (1.4.4) is well-known in the $X$ smooth case.) Recall that we have by (1.1.3)
$$R_y(\beta_j)=\frac{\beta_j}{Q_y(\beta_j)}.$$
We thus get
$$i_*\bigl(T_{y*}^{\vir}(X)\bigr)=\bigl(\mprod_j\,R_y(\beta_j)\cup T^*_y(TY)\bigr)\cap[Y].
\leqno(1.4.5)$$
By applying the trace morphism ${\rm Tr_Y}:H_{\ssb}(Y)\to\Q$, this implies the compatibility with Hirzebruch's construction \cite[11.2]{Hi} 
$${\rm Tr}_X\bigl(T_{y*}^{\vir}(X)\bigr)=\int_Y\mprod_j\,R_y(\beta_j)\cup T^*_y(TY).
\leqno(1.4.6)$$
Here $\int_Y:H^{2\dim Y}(Y)\to\Q$ denotes the canonical morphism (which is also called the trace morphism), and the right-hand side of (1.4.6) is equal to $T_y(\beta_1,\dots,\beta_r)_Y$ in the notation of \cite[11.2]{Hi} where $Y$ and $\beta_j$ are respectively denoted by $M$ and $v_j$.
\sk
Note that the above argument is mostly useful for the degree zero part of the (homology) Hirzebruch class unless the natural morphism $i_*:\HH_{\ssb}(X)\to\HH_{\ssb}(Y)$ is injective since the information may be lost in the other case.
\ms\nin
{\bf 1.5.~Spectrum.} Let $f$ be a holomorphic function on a complex manifold $Y$ of dimension $n$. Let $x\in X:=f^{-1}(0)\subset Y$. We have the Steenbrink spectrum
$$\aligned\Sp(f,x)&=\msum_{\al\in\Q}\,n_{f,x,\al}\,t^{\al}\q\h{with}\q n_{f,x,\al}:=\msum_j\,(-1)^{j-n+1}h_{f,x,\la}^{p,j-p},\\h_{f,x,\la}^{p,j-p}:&=\dim\Gr^p_F\Ht^j(F_{f,x},\C)_{\la}\q\bl(p=\lfloor n-\al\rfloor,\q\la=\exp(-2\pi i\al)\br).\endaligned
\leqno(1.5.1)$$
Here $F$ is the Hodge filtration of the canonical mixed Hodge structure on the reduced Milnor cohomology $\Ht^j(F_{f,x},\C)$ with $F_{f,x}$ the Milnor fiber of $f$ around $x$, and $\Ht^j(F_{f,x},\C)_{\la}$ is the $\la$-eigenspace of the cohomology by the semisimple part $T_s$ of the Milnor monodromy $T$, see \cite{St1}, \cite{St2}. Recall that
$$\lfloor\al\rfloor:=\max\{i\in\Z\mid i\les\al\},\q\lceil\al\rceil:=\min\{i\in\Z\mid i\ges\al\}.
\leqno(1.5.2)$$
\sk
Let $i_x:\{x\}\into X$ denote the inclusion. Then we have an isomorphism of mixed Hodge structures
$$\Ht^j(F_{f,x},\Q)=H^ji_x^*\varphi_f\Q_{h,Y},
\leqno(1.5.3)$$
which is compatible with the action of the semisimple part $T_s$ of the Milnor monodromy $T$. Here the category of mixed Hodge modules on a point is identified with the category of graded-polarizable mixed $\Q$-Hodge structures \cite{De2}, see \cite{Sa2}. In fact, (1.5.3) is actually the definition of the mixed Hodge structure on the left-hand side.
\bs\bs
\vbox{\centerline{\bf 2. Proofs of the main assertions}
\bs\nin
In this section we give the proofs of Theorem~1 and Propositions~1 and 2.}
\ms\nin
{\bf 2.1.~Proof of Proposition~1.} Fix a base point $s_0\in S$. Associated with the local system $\Hc^j_{S,\la}$, we have the monodromy representation
$$\rho^j_{S,\la}:\pi_1(S,s_0)\to{\rm Aut}(\Hc^j_{S,\la,s_0}).
\leqno(2.1.1)$$
Any $\gamma\in\pi_1(S,s_0)$ is represented by a piecewise linear path (using local coordinates). It is contractible inside $\St$ by condition~$(c)$. We may assume that this contraction is also given by a piecewise linear one, and intersects $D_{\St}=\St\setminus S$ transversally at smooth points. (Here it may be better to use a sufficiently fine triangulation of $\St$ compatible with $D_{\St}$.) Then, using condition~$(d)$, we see that the image of $\gamma$ by $\rho^j_{S,\la}$ is given by the multiplication by
$$\mprod_i\,(c_{S,i,\la})^{a_i(\gamma)}\in\C^*,
\leqno(2.1.2)$$
since the multiplication by $c_{S,i,\la}$ belongs to the center of ${\rm Aut}(\Hc^j_{S,\la,s_0})$. Here $a_i(\gamma)\in\Z$ depends only on the contraction of $\gamma$, and is independent of $j$.
\sk
This argument implies that we have a monodromy representation
$$\rho_{S,\la}:\pi_1(S,s_0)\to\C^*\,\bl(={\rm Aut}(\C)\br),
\leqno(2.1.3)$$
such that any nonzero $\rho^j_{S,\la}$ is isomorphic to a direct sum of copies of $\rho_{S,\la}$.
Let $L_{S,\la}$ be the local system corresponding to $\rho_{S,\la}$. Then the assertion follows. This finishes the proof of Proposition~1.
\ms\nin
{\bf 2.2.~Proof of Proposition~2.} By condition~(0.4) together with (2.1.2), the monodromy representation (2.1.3) is compatible with the product between the $\la\in\La_S$. Note that $\La_S$ is stable by inverse. (In fact, the local systems $\Hc^j_S$ are defined over $\Q$ so that $\La_S$ is stable by complex conjugation, and the eigenvalues of the Milnor monodromies are roots of unity.) We then get the monodromy representation for any $\la\in G_S$
$$\rho_{S,\la}:\pi_1(S,s_0)\to\C^*,
\leqno(2.2.1)$$
in a compatible way with the product between the $\la\in G_S$.
We define $L'_S$ to be the rank 1 local system corresponding to $\rho_{S,\e(1/m'_S)}$.
\sk
By conditions~(0.4) and (0.6) the eigenvalues of the residues of the connection of $\Lc_{\St}^{\prime\,\otimes k}$ along $D_{{\St},i}$ is given by
$$km'_{S,i}/m'_S.
\leqno(2.2.2)$$
In fact, we have $m'_{S,i}/m'_S\in[0,1)$ for $k=1$ by condition~(0.6). So (0.8) follows, and the assertion~(i) is proved. The argument is similar for the assertion~(ii). This finishes the proof of Proposition~2.
\ms\nin
{\bf 2.3.~Proof of Theorem~1.}
Let $\M^j_{\St,\la}$ be the Deligne extension of the local system $\Hc^j_{S,\la}$ such that the eigenvalues of the residues of the connection are contained in $(0,1]$ (see \cite{De1}).
By \cite{MSS}, Propositions~5.1.1 and 5.2.1, we have
$$\aligned&\q T_{y*}\bl((i_{\Si\setminus X',\Si})_{!\,}\varphi_f\Q_{h,Y\setminus X'}\br)\\
&=\sum_{S,j,p,q,\la}(-1)^{j+q}\,(\pi_{\St})_*\,td_{(1+y)*}\bl[\Gr_F^p\M^j_{\St,\la}\otimes_{\Oc_{\St}}\Omega^q_{\St}(\log D_{\St})\br](-y)^{p+q},\endaligned
\leqno(2.3.1)$$
By conditions~$(a)$ and $(b)$ together with (1.5.1) and (1.5.3), we then get
$$\aligned&\q T_{y*}\bl((i_{\Si\setminus X',\Si})_{!\,}\varphi_f\Q_{h,Y\setminus X'}\br)\\
&=\sum_{S,j,p,q,\la}(-1)^{j+q}\,h_{f,S,\la}^{p,j-p}\,(\pi_{\St})_*\,td_{(1+y)*}\bl[\Lc_{\St,\la}\otimes_{\Oc_{\St}}\Omega^q_{\St}(\log D_{\St})\br](-y)^{p+q},\endaligned
\leqno(2.3.2)$$
with
$$h_{f,S,\la}^{p,j-p}:=h_{f,x,\la}^{p,j-p}\q\h{for any}\,\,\,x\in S.$$
\sk
In fact, condition~$(a)$ implies that the Hodge filtration $F$ is defined on the level of local systems, and the graded pieces of the filtration $F$ are still direct sums of rank 1 local systems as in condition~$(b)$. The assertion now follows from (2.3.2) and (1.5.1). This finishes the proof of Theorem~1.
\ms\nin
{\bf 2.4.~Proof of \cite[Proposition~4.1]{MSS} in the case $m>1$.}
By the same argument as in \cite{MSS}, the assertion is reduced to the normal crossing case. Then it is enough to show the vanishing of the reduced cohomology of
$$U_{\ep,t}:=\bl\{(y_1,\dots,y_n)\in\C^n\,\big|\,\msum_{i=1}^{n-1}|y_i|^2<\ep^2,\,|y_n|<\ep,\,g=y_n^m\,t\br\}\q(0<|t|\ll\ep\ll 1),$$
by using the fundamental neighborhood system of $0\in\C^n$ given by
$$U_{\ep}:=\bl\{(y_1,\dots,y_n)\in\C^n\,\big|\,\msum_{i=1}^{n-1}|y_i|^2<\ep^2,\,|y_n|<\ep\br\}\q(0<\ep\ll 1),$$
where $g=\prod_{i=1}^ry_i^{m_i}$ with $r<n$. Consider
$$V_{\ep,t}:=\bl\{(y_1,\dots,y_{n-1})\in\C^{n-1}\,\big|\,\,\msum_{i=1}^{n-1}|y_i|^2<\ep^2,\,|g|<\ep^m|t|\br\}\q(0<|t|\ll\ep\ll 1).$$
It is contractible (see \cite{Mi}), and $U_{\ep,t}$ is a ramified covering of $V_{\ep,t}$ which is ramified over the normal crossing divisor $V_{\ep,t}\cap g^{-1}(0)$. Then $V_{\ep,t}$ and $U_{\ep,t}$ retract to $V_{\ep,t}\cap g^{-1}(0)$, and the assertion follows.
\ms\nin
{\bf 2.5.~Formula for the Chern-Milnor classes.}
Let $M(X)\in H_{\ssb}(\Si)$ be the Chern-Milnor class $M(X)$ as in \cite{PP} (which can be obtained by specializing $M_y(X)$ to $y=-1$).
Then, without assuming conditions~$(a)$, $(b)$, we have
$$\aligned M(X)&=\msum_{S\in\Sc}\,\widetilde{\chi}(F_{f,S})\,c_*(1_S)\\&=\msum_{S\in\Sc}\,\widetilde{\chi}(F_{f,S})\,(\pi_{\St})_*\bl(c^*\bl(\Omega_{\St}^1(\log D_{\St})^{\vee}\br)\cap[\St]\br).\endaligned
\leqno(2.5.1)$$
Here $c_*(1_S)\in H_{\ssb}(\Si)$ is as in \cite{Ma}, and
$$\widetilde{\chi}(F_{f,S}):={\chi}(F_{f,S})-1,$$
with $F_{f,S}$ the Milnor fiber of $f$ around a sufficiently general $x\in S$. The second equality of (2.5.1) follows from \cite{Al}, \cite{GP}.
Note that the $\widetilde{\chi}(F_{f,S})$ ($S\in\Sc$) give the constructible function associated with the vanishing cycle complex $\varphi_f\Q_{Y\setminus X'}$.
\sk
It is well-known that, if the restriction of $X$ to a transversal slice to $S$ at $x$ is locally defined by a homogeneous polynomial $g_S$ of degree $m_S$, then
$${\chi}(F_{f,S})=\chi\bl(\PP^{\,c_S-1}\setminus g_S^{-1}(0)\br)\,m_S,
\leqno(2.5.2)$$
where $c_S:=\codim\,S$.
\sk
In the hyperplane arrangement case, it is known (see \cite{STV}) that
$$\h{$\chi\bl(\PP^{\,c_S-1}\setminus g_S^{-1}(0)\br)=0$ if and only if $\So$ is not a dense edge,}$$
in the notation of (3.1) below.
\bs\bs
\vbox{\centerline{\bf 3. Hyperplane arrangement case}
\bs\nin
In this section we treat the hyperplane arrangement case.}
\ms\nin
{\bf Notation~3.1.}
Assume $X$ is a projective hyperplane arrangement in $Y=\PP^n$, where $L=\Oc_{\PP^n}$ and $X'$ is a sufficiently general hyperplane. Let $X_j$ be the irreducible components of $X$ with multiplicities $m_j$ ($j=1,\dots,r$). Note that
$$m=\msum_j\,m_j,$$
and
$$X_j\subset\Si\q\h{if}\q m_j>1.$$
\sk
We have a canonical stratification $\Sc=\Sc_{\Si'}$ of $\Si':=\Si\setminus X'$ such that
$$\So=\mcap_{j\in I(S)}\,X_j,\q S=\So\setminus\bl(\mcup_{j\notin I(S)}\,X_j\cup X'\br)\q\h{with}\q I(S):=\{j\mid X_j\supset S\}.
\leqno(3.1.1)$$
For the proof of Proposition~3, we have to consider also the canonical stratification $\Sc_X$ of $X$ such that (3.1.1) holds by deleting $X'$.
Here $\So$ is called an {\it edge} of the hyperplane arrangement $X$.
\sk
Let $C(X)$ denote the corresponding central hyperplane arrangement of $V:=\C^{n+1}$ with irreducible components $C(X_j)$ and multiplicities $m_j$. Here $C(X_j)$ denotes the cone of $X_j$. Set
$$V^S:=V/V_S\q\h{with}\q V_S:=C(\So).$$
\sk
For each $S\in\Sc$, we have the quotient central hyperplane arrangement
$$C(X)^S\subset V^S$$
defined by the affine hyperplanes
$$C(X_j)^S:=C(X_j)/V_S\subset V^S\q\h{for}\,\,\,j\in I(S),$$
where the irreducible components $C(X_j)^S$ have the induced multiplicities $m_j$.
\sk
For each $S\in\Sc$, we also have the induced projective hyperplane arrangement
$$X_{\So}:=\mcup_{j\notin I(S)}\,X_j\cap\So\subset\So,$$
such that $X_{\So}\setminus\Si$ has the induced stratification
$$\Sc_S:=\{S'\in\Sc\mid S'\subset X_{\So}\}.$$ 
Here each $S'\in\Sc_S$ (especially for $\codim_{\So}S'=1$) has the induced multiplicity
$$m_{S',S}:=\msum_{j\in I(S')\setminus I(S)}\,m_j.$$
\sk
Let $f^S$ be a homogeneous polynomial defining $C(X)^S\subset V^S$. This is essentially identified with $g_S$ in (0.5), and
$$m_S=\deg f^S=\msum_{j\in I(S)}\,m_j\in\Z\,m'_S,
\leqno(3.1.2)$$
where $m'_S$ is as in (0.6).
Note that there is a shift of indices of the spectral numbers
$$n_{f,S,\al}=(-1)^{\dim S}n_{f^S,0,\beta}\q\h{with}\,\,\,\beta=\al-\dim S,
\leqno(3.1.3)$$
and a formula for $n_{f^S,0,\beta}$ in the reduced case with $\dim V^S\les 3$ can be found in \cite{BS2}.
\sk
As for the good compactification $\St$ of $S\in\Sc$, it can be obtained by blowing-up $\So\subset\PP^n$ along the edges $\So'$ of $X_{\So}\subset\So$ with $\codim_{\So}\So'\ges 2$ by decreasing induction on the codimension of the edges, where we restrict to the $S'$ such that $X_{\So}$ is not a divisor with normal crossings on any neighborhood of $S'$ in $\So$ as in \cite{BS2}. (However, we do not restrict to the dense edges as in \cite{STV}, see Remark~(3.4)(i) below for the definition of dense edge.)
\sk
Let $\E_{\So',\St}$ be the proper transform of the exceptional divisor of the blow-up along $\So'$ in $\So$ if $\codim_{\So}\So'\ges 2$ (where we assume that $X_{\So}\subset\So$ is not a divisor with normal crossings on any neighborhood of $S'$ as above).
If $\codim_{\So}\So'=1$, then $\E_{\So',\St}$ denotes the proper transform of $\So'\subset\So$.
Let $\X_{\infty,\St}$ be the proper transform of $X_{\infty,\So}:=X'\cap\So\subset\So$. These are the irreducible components $D_{{\St},i}$ of $D_{\St}\subset\St$. The integers $m'_{S,i}$ in (0.4) for the components $\E_{\So',\St}$ and $\X_{\infty,\St}$ will be denoted respectively by
$$m'_{S',S},\q m'_{\infty,S}.$$
\ms\nin
{\bf Proposition~3.2.} {\it With the above notation, condition~$(0.4)$ holds with
$$m'_{S',S}/m_S=\{m_{S',S}/m_S\},\q m'_{\infty,S}/m_S=\{-m/m_S\},
\leqno(3.2.1)$$
where $\{\al\}:=\al-\lfloor\al\rfloor$. Moreover, $L_S$ in Proposition~{\rm 2\,(ii)} exists, and we have
$$\Lc_{\St}=\pi_{\St}^*\,\Oc_{\So}\bl(-\bl\lceil\,\msum_{j\notin I(S)}\,m_j/m_S\,\br\rceil\br)\otimes_{\Oc_{\St}}\Oc_{\St}\bl(\msum_{S'\subset S}\,\lfloor m_{S',S}/m_S\rfloor\E_{\So',\St}\br),
\leqno(3.2.2)$$
where $\Oc_{\So}(k)$ denotes the pull-back of $\Oc_{\PP^n}(k)$ by $\So\into\PP^n$.}
\ms\nin
{\it Proof.} By using the blowing-ups along the $\So'$ in $\PP^n$, the assertion (3.2.1) is reduced to Lemma~(3.3) below.
The integrable connection corresponding to the local system $L_S$ can be constructed easily on $S\subset\So\setminus X'\cong\C^{\dim S}$, see e.g. \cite{ESV}.
\sk
For the proof of (3.2.2), it is then enough to show
$$\Lc_{\St}^{\otimes m_S}=\pi_{\St}^*\,\Oc_{\So}\bl(-m_S\bl\lceil\,\msum_{j\notin I(S)}\,m_j/m_S\,\br\rceil\br)\otimes_{\Oc_{\St}}\Oc_{\St}\bl(\msum_{S'\subset S}\,m_S\lfloor m_{S',S}/m_S\rfloor\E_{\So',\St}\br).
\leqno(3.2.3)$$
(Indeed, the simply connectedness implies that $\CH^1(\St)$ is torsion-free, since it implies that $H^1(\St,\Z)=0$ and $H^2(\St,\Z)$ is torsion-free.)
The left-hand side of (3.2.3) is a line bundle with a logarithmic connection such that the eigenvalues of the residues along $\E_{\So',\St}$ and $\X_{\infty,\St}$ are respectively $m'_{S',S}$ and $m'_{\infty,S}$ (since $m'_{S',S}/m_S,\,m'_{\infty,S}/m_S\in[0,1)$ by (3.2.1)).
So we get
$$\Lc_{\St}^{\otimes m_S}=\Oc_{\St}\bl(-\msum_{S'\subset S}\,m'_{S',S}\E_{\So',\St}-m'_{\infty,S}\X_{\infty,\St}\br).$$
On the other hand, setting $X_{j,\So}:=X_j\cap\So$, we have
$$\aligned&\pi_{\St}^*\,\Oc_{\So}\bl(-m_S\bl\lceil\,\msum_{j\notin I(S)}\,m_j/m_S\,\br\rceil\br)\\&=\pi_{\St}^*\,\Oc_{\So}\bl(-\msum_{j\notin I(S)}\,m_jX_{j,\So}-m'_{\infty,S}X_{\infty,\So}\br)\\&=\Oc_{\St}\bl(-\msum_{S'\subset S}\,m_{S',S}\E_{\So',\St}-m'_{\infty,S}\X_{\infty,\St}\br).\endaligned$$
Here the first isomorphism follows from the second equality of (3.2.1) (which implies that $m'_{\infty,S}\in[0,m_S)$), since we have by the definition of $m$ and $m_S$
$$\msum_{j\notin I(S)}\,m_j=m\,\,\,\,\h{mod}\,\,\,\,m_S.$$
The second isomorphism is obtained by calculating the total transform of the divisor. So the assertion follows from (3.2.1) which implies that
$$m'_{S',S}=m_{S',S}-m_S\lfloor m_{S',S}/m_S\rfloor.$$
\ms\nin
{\bf Lemma~3.3.} {\it Let $g$ be a holomorphic function on a complex manifold $Y$. Set $X:=Y\times\C^*$, and $f:=gz^a$ where $z$ is the coordinate of $\C$ and $a\in\Z$. Then the monodromy of the local system $\Hc^j\psi_f\Q_X|_{\{y\}\times\C^*}$ for $y\in Y$ is given by $T^{-a}$, where $T$ is the Milnor monodromy.}
\ms\nin
{\it Proof.} Since $f^{-1}(0)$ is analytic-locally trivial along $\{y\}\times\C^*$, we can calculate $\Hc^j(\psi_f\Q_X)_{(y,z)}$ by the cohomology of
$$\bl\{\,y'\in Y\,\,\big|\,\,||y'||<\ep,\,g(y')=z^{-a}t\,\br\}\q(0<|t|\ll\ep\ll 1).$$
Here $||y'||$ is defined by taking local coordinates of $Y$ around $y$, and we may assume $|z|=1$ for the calculation of the monodromy. Then the assertion is clear. This finishes the proofs of Lemma~(3.3) and Proposition~(3.2).
\ms\nin
{\bf Remarks~3.4.} (i) An edge $\So$ of a hyperplane arrangement $X$ is called {\it dense} (see \cite{STV}) if its associated quotient central hyperplane arrangement $C(X)^S$ is indecomposable. Note that a central hyperplane arrangement is called indecomposable if it is not a union of hyperplane arrangements coming from $\C^{n_1}$ and $\C^{n_1}$ via the projections
$$\C^n\to\C^{n_1},\q\C^n\to\C^{n_2},$$
where $n=n_1+n_2$ and $n_1,n_2>0$.
\ms
(ii) If $Y$ is a simply connected smooth variety and $Z\subset Y$ is a closed subvariety of codimension at least two, then $Y\setminus Z$ is simply connected. (Indeed, a contraction of a path has real dimension 2, and can be modified so that it does not intersect $Z$.)
\sk
This implies that if $Y$ is a simply connected smooth variety and $Y'\to Y$ is a proper birational morphism from a smooth variety, then $Y'$ is also simply connected. (Indeed, if a smooth variety has a dense Zariski-open subvariety which is simply connected, then it is also simply connected. This follows from the fact that a path has real dimension 1, and may be modified so that it is contained in the open subvariety.)
\ms
(iii) The image of the cycle map is contained in
$$\Gr^W_{-2k}H_{2k}(\Si,\Q)\subset H_{2k}(\Si,\Q),$$
and so is the image of $td_*$. Moreover, the structures of $\Gr^W_{-2k}H_{2k}(\Si,\Q)$ and $\CH_k(X)_{\Q}$ in the projective hyperplane arrangement case are quite simple as below.
\ms\nin
{\bf Proposition~3.5.} {\it Let $X$ be a projective hyperplane arrangement in $\PP^n$ with $X_j\,\,(j\in[1,r])$ the irreducible components. Let $W$ be the weight filtration of the canonical mixed Hodge structure on $H_k(X,\Q)$ $($which is the dual of $H^k(X,\Q))$. Then we have
$$\Gr^W_{-k}H_k(X,\Q)=\begin{cases}\mopl_{1\les j\les r}\,\Q[X_j]&\h{if}\,\,\,k=2n-2,\\ \Q(k/2)&\h{if}\,\,\,k\in 2\,\Z\cap[0,2n-4],\\ \,0&\h{otherwise,}\end{cases}
\leqno(3.5.1)$$
where $[X_j]$ denotes the class of $X_j\,\,($and $\Q[X_j]$ is not a polynomial algebra$)$. Moreover we have the canonical isomorphisms}
$$\Gr^W_{-k}H_k(X,\Q)\simto H_k(\PP^n,\Q)(k/2)\q(0\les k<2n-2).
\leqno(3.5.2)$$
\ms\nin
{\it Proof.} We have a long exact sequence
$$\to H^{\BM}_{k+1}(\PP^n\setminus X,\Q)\to H_k(X,\Q)\to H_k(\PP^n,\Q)\to,$$
together with
$$H^{\BM}_{k+1}(\PP^n\setminus X,\Q)=H^{2n-k-1}(\PP^n\setminus X,\Q)(n).$$
Here $H^p(\PP^n\setminus X,\Q)$ vanishes for $p>n$, and has type $(p,p)$ for $p\les n$ by \cite{Br}. (In fact, an integral logarithmic $p$-form has type $(p,p)$ by \cite[Thm.~8.2.4\,(i)]{De4}, since the latter implies that it induces a morphism of mixed Hodge structures from $\Q(-p)$ to the $p$-th cohomology group.)
\sk
Setting $p=2n-k-1$, we see that $H^{\BM}_{k+1}(\PP^n\setminus X,\Q)=0$ if $k+1<n$, and it has type
$$(n-k-1,n-k-1),\,\,\,\h{if}\,\,\,\,k+1\ges n.$$
Here we have
$$2n-2k-2>-k\q\h{for}\q k\in[0,2n-3].$$
So we get (3.5.2) and also (3.5.1) except for $k=2n-2$.
In the last case, the above argument shows that $H_{2n-2}(X,\Q)$ has type $(1-n,1-n)$, and has dimension $r$ by using \cite{Br}. So the assertion follows.
\ms\nin
{\bf Proposition~3.6.} {\it With the notation of $(3.1)$ and Proposition~$(3.5)$ above, set
$$\aligned\Si_1&:=\mcup_{X_j\subset\Si}\,X_j,\q\Si_2:=\overline{\Si\setminus\Si_1},\q\Sc^{(i)}:=\{S\in\Sc\mid\codim_YS=i\},\\
\Sc_a&:=\{S\in\Sc\mid S\in\Si_a\},\q\Sc_a^{(i)}:=\Sc_a\cap\Sc^{(i)}\,\,\,\,(a=1,2).
\endaligned$$
Then we have
$$\CH_k(X)_{\Q}=\begin{cases}\mopl_{1\les j\les r}\,\Q[X_j]&\h{if}\,\,\,k=n-1,\\ \Q&\h{if}\,\,\,k\in[0,n-2].\end{cases}
\leqno(3.6.1)$$
$$\CH_k(\Si)_{\Q}=\begin{cases}\mopl_{S\in\Sc_1^{(1)}}\,\Q[\So]&\h{if}\,\,\,k=n-1,\\ \mopl_{S\in\Sc_2^{(2)}}\,\Q[\So]\oplus\Q&\h{if}\,\,\,k=n-2\,\,\,\h{and}\,\,\,\Si_1\ne\emptyset,\\ \mopl_{S\in\Sc_2^{(2)}}\,\Q[\So]&\h{if}\,\,\,k=n-2\,\,\,\h{and}\,\,\,\Si_1=\emptyset,\\ \Q&\h{if}\,\,\,k\in[0,n-3].\end{cases}
\leqno(3.6.2)$$
Moreover, for $S\in\Sc^{(i)}$ with $i\ges 3$, we have the canonical isomorphisms
$$\CH_k(\So)_{\Q}\simto\CH_k(\Si)_{\Q}\simto\CH_k(\PP^n)_{\Q}\q(0\les k\les\dim\So\les n-3),
\leqno(3.6.3)$$
$$\CH_{n-2}(\So)_{\Q}\simto\CH_{n-2}(\Si_1)_{\Q}\simto\CH_k(\PP^n)_{\Q}\q\h{if}\q S\in\Sc^{(2)}_1,
\leqno(3.6.4)$$
and $\Q$ in $(3.6.2)$ for $k=n-2$ and $\Si_1\ne\emptyset$ is given by the image of $\CH_{n-2}(\Si_1)_{\Q}$ in $(3.6.4)$.}
\ms\nin
{\it Proof.} Since $\dim\Si_1=n-1$ and $\dim\Si_2=n-2$, the assertions easily follow from the well-known facts that we have for $\So'\subset\So\subset\PP^n$
$$\CH_k(\So)=\Q\q\h{for}\,\,\,k\in[0,\dim\So],
\leqno(3.6.5)$$
$$\CH_k(\So')\simto\CH_k(\So))\simto\CH_k(\PP^n)\,\,\,\h{for}\,\,\,k\in[0,\dim\So'].
\leqno(3.6.6)$$
For instance, (3.6.6) implies that the image of $[\So']$ in $\CH_{n-2}(\Si_1)_{\Q}$ is independent of $\So'\in\Sc_1^{(2)}$ (by applying it to the $\So\in\Sc_1^{(1)}$). So (3.6.4) follows. The proofs of the other assertions are similar. This finishes the proof of Proposition~(3.6).
\ms\nin
{\bf Proposition~3.7.} {\it The Hirzebruch-Milnor class $M_y(X)$ of an hyperplane arrangement $X$ in $\PP^n$ is a combinatorial invariant, where $\HH_{\ssb}(\Si)=\CH_{\ssb}(\Si)_{\Q}$ and Proposition~$(3.6)$ is used.}
\ms\nin
{\it Proof.} Let $\Ec_{\St,\al,q}$ be a vector bundle on $\St$ defined by
$$\Ec_{\St,\al,q}:=\Lc_{\St,\e(-\al)}\otimes_{\Oc_{\St}}\Omega^q_{\St}(\log D_{\St}).$$
By Theorem~1, Proposition~(3.6) and \cite{BS2} (which implies that the $n_{f,S,\al}$ are combinatorial invariants) together with the definition of $td_{(1+y)*}$ in (1.2.1), it is enough to show that the following is a combinatorial invariant for any $S,\al,q:$
$$(\pi_{\St,\So})_*td_*(\Ec_{\St,\al,q})\in\HH_{\ssb}(\So)=\Q^{\dim\So+1}.\leqno(3.7.1)$$
Here $\pi_{\St,\So}:\St\to\So$ is the canonical morphism, and (3.6.5) is used for the last isomorphism.
Since $\So$ is projective space, we may assume for the proof of (3.7.1)
$$\HH_{\ssb}(\So)=H_{2\ssb}(\So,\Q),\q\HH^{\ssb}(\So)=H^{2\ssb}(\So,\Q),$$
and similarly for $\St$. By \cite{BFM} we have
$$td_*(\Ec_{\St,\al,q})=\bl(ch(\Ec_{\St,\al,q})\cup td^*(T\St)\br)\cap[\St]\q\h{in}\q\HH_{\ssb}(\St).$$
This means that $td_*(\Ec_{\St,\al,q})$ is identified by Poincar\'e duality with
$$ch(\Ec_{\St,\al,q})\cup td^*(T\St)\in\HH^{\ssb}(\St).$$
Moreover, the pushforward by $\pi_{\St,\So}$ is calculated by the top degree part of
$$ch(\Ec_{\St,\al,q})\cup td^*(T\St)\cup\pi_{\St,\So}^*e^k,$$
where $e^k$ is the canonical generator of $\HH^k(\So)$ ($k\in[0,\dim\So]$). These can be calculated by the method in \cite{BS2}, Section 5 using the combinatorial description of the cohomology ring as in \cite{DP}.
Here $[\E_{\So',\St}]$, $[\pi_{\St,\So}^*e]\in\HH^1(\St)$ correspond to $e_V\,\,(V\ne 0)$ and $-e_0$ in \cite[5.3]{BS2}. Moreover $c^*(\Omega_{\St}^q(\log D_{\St})$ and $td^*(T\St)$ are combinatorially expressed (loc.~cit., 5.4) by using certain universal polynomials.
Combining these with Propositions~2\,(ii) and (3.2), it follows that (3.7.1) is a combinatorial invariant. Here it is not necessary to use the full result of \cite{DP}, since we need only the fact that the multiple-intersection numbers of the $e_V$ are independent of the position of the irreducible components, and are determined combinatorially (together with certain relations among the $e_V$). This finishes the proof of Proposition~(3.7).

\end{document}